\documentclass[12pt]{article}
\usepackage{amssymb, amsmath}
\begin{document}
\voffset=0.0truein \hoffset=-0.5truein
\setlength{\textwidth}{6.0in} \setlength{\textheight}{8.8in}
\setlength{\topmargin}{-0.2in}
\renewcommand{\theequation}{\arabic{section}.\arabic{equation}}
\newtheorem{theorem}{Theorem}[section]
\newtheorem{lemma}{Lemma}[section]
\newtheorem{pro}{Proposition}[section]
\newtheorem{cor}{Corollary}[section]
\newcommand{\n}{\nonumber}
\newcommand{\intr}{\int_{\Bbb R^2}}
\newcommand{\x}{{X_\alpha}}
\newcommand{\y}{{Y_\alpha}}
\newcommand{\lo}{\lambda_1}
\newcommand{\lt}{\lambda_2}
\newcommand{\ltt}{\lambda_3}
\newcommand{\lf}{\lambda_4}
\newcommand{\lfr}{\frac{\lambda_4}{\lambda_2}}
\newcommand{\lfrr}{\frac{2\lambda_4}{\lambda_2}}
\renewcommand{\do}{\delta_1}
\newcommand{\dt}{\delta_2}
\renewcommand{\a}{\alpha}
\newcommand{\ko}{\kappa_1}
\newcommand{\kt}{\kappa_2}
\newcommand{\aff}{{a\l}^{\frac{1}{a}}}
\newcommand{\af}{^{\frac{1}{a}}}
\newcommand{\vare}{\varepsilon}
\newcommand{\ro}{\rho^{I}_{ a, \varepsilon} (z)}
\newcommand{\rt}{\rho^{II}_{a, \varepsilon} (z)}
\newcommand{\go}{g^{I}_{a,\varepsilon}(z)}
\newcommand{\gt}{g^{II}_{a,\varepsilon }(z)}
\newcommand{\rro}{\rho_1 }
\newcommand{\rrt}{\rho_2 }
\newcommand{\rn}{r^{2N+2}}
\newcommand{\rr}{r^{2}}
\newcommand{\orn}{(1+r^{2N+2})}
\newcommand{\orr}{(1+r^{2})}
\newcommand{\intt}{\int_0 ^\infty }
\newcommand{\php}{\varphi_{+}}
\newcommand{\phm}{\varphi_{-}}
\newcommand{\bo}{\beta_1}
\newcommand{\bt}{\beta_2}
\newcommand{\ub}{\hat{u}}
\newcommand{\eb}{\hat{\eta}}
\newcommand{\ut}{\tilde{u}}
\newcommand{\et}{\tilde{\eta}}
\renewcommand{\n}{\nonumber}
\newcommand{\bb}{\begin{equation}}
\newcommand{\ee}{\end{equation}}
\newcommand{\bq}{\begin{eqnarray}}
\newcommand{\eq}{\end{eqnarray}}
\newcommand{\bqn}{\begin{eqnarray*}}
\newcommand{\eqn}{\end{eqnarray*}}
\title{Existence of Multistring Solutions of the Self-Gravitating Massive $W-$Boson}
\author{Dongho Chae\thanks{Permanent Address: Department of Mathematics,
              Sungkyunkwan University, Suwon 440-746, Korea,
              e-mail : chae@skku.edu}\\
Center for Scientific Computation And\\
 Mathematical Modeling\\
 Paint Branch Drive \\
 University of Maryland\\
  College Park, MD 20742-3289, USA\\
  e-mail: {\it dchae@cscamm.umd.edu}}
\date{}
\maketitle
\begin{abstract}
 We consider a semilinear elliptic system which include the model
 system of the $W-$strings in the cosmology as a special case.  We
 prove existence of multi-string solutions and obtain precise
 asymptotic decay estimates near infinity for the solutions.
 As a special case of this result we solve an open problem posed in
 \cite{yan}
 \\
 \ \\
{\textbf{Key Words:} semilinear elliptic system,  exponential
nonlinearities, selfdual gauge field theories}\\
{\textbf{AMS Subject Classification:} 35J45, 35J60, 37K40, 70S15.}
\end{abstract}
\section{Introduction}
 \setcounter{equation}{0}
 Let $\lo ,\lt ,\lt, \lf  >0$ be given.
 We consider the following system for $(u, \eta )$ in $\Bbb R^2$.
 \bq
 \label{11}
 \Delta u&=&-\lo e^{\eta} -\lt e^{u}+4\pi \sum_{j=1} ^{N} \delta (z-z_{j} ),\\
 \label{12}
 \Delta \eta &=& -\ltt e^{\eta} -\lf e^{u}
 \eq
 equipped with the boundary condition
 \bb
 \label{13}
 \intr e^{u}dx +\intr e^{\eta} dx < \infty ,
 \ee
  where we denoted $z=x_1+ix_2 \in \Bbb C=\Bbb R^2$.
  The system
  (\ref{11})-(\ref{12}) is the reduced form of the  Bogomol'nyi type of
  equations
   modeling the cosmic strings with matter field given by the
   massive $W-$boson of the electroweak theory, if we choose the coefficients as,
   \bb
   \label{13a}
   \lo=2m_W^2 , \quad \lt=4e^2, \ltt =\frac{16\pi G
   m_W^4}{e^2},\quad \lf=32\pi G m_W ^2,
   \ee
   where $m_W$ is the mass of the $W-$boson, $e$ is the charge of
   the electron, and $G$ is the gravitational constant.
   The points $\{ z_1, \cdots,
   z_N\}$ corresponds to the location on the $(x_1,x_2)-$plane
   of parallel (along the $x_3-$axis) strings. See \cite{yan,amb}
   for the derivation of this system from the corresponding
   Einstein-Weinberg-Salam theory as well as interesting physical backgrounds of the model.
   There are many previous mathematical studies on the planar
   electroweak theory recently(\cite{spr1,spr2,bar2,
   cha2}). In particular in \cite{cha3} the authors considered full
   electroweak field as the matter field coupled with the
   gravitation. In the model from which our system is derive the
   matter field coupled to gravity is the massive $W-$boson.
    In \cite{yan} the construction of radially symmetric solutions(in
    the case $z_1=\cdots=z_N$) of (\ref{11})-(\ref{13}) is discussed by
    further reduction
    the system into a single equation, and solving the ordinary
    differential equation.
    When the locations of strings are different to each other, however, we
    cannot assume the radial symmetry of the solutions,
    and no existence theory is available. In particular, the
    author of \cite{yan} left the  construction of solution  in this case as an open
   problem. One of our main purpose in
   this paper is to solve this problem. Actually, we solve the
   existence problem for more general coefficient cases as in
   (\ref{11})-(\ref{12}). The following is our main theorem.
  \begin{theorem}
 Let $N \in \Bbb N \cup \{0\}$, and $ \mathcal{Z}=\{ z_{j}\}_{j=1} ^{N} $
 be given in $\Bbb R^2$
 allowing multiplicities. Suppose the coefficients, $\lambda_1,
 \lambda_2,\lambda_3,\lambda_4 $ satisfy  one of the conditions;
either
 \bb\label{condition}
\lo\lf -\lt\ltt=0, \ee
 or
 \bb\label{conditionA}
 \lo\lf -\lt\ltt\neq 0 \quad \mbox{and}\quad \frac{ \lt}{2\lf} < N+1.
 \ee
  Then, there exists a
constant $\vare_1
>0$
 such that
 for any $\vare \in (0, \vare_1 )$ and any $c_0 >$ there exists a family of  solutions to
 (\ref{11})-(\ref{13}),
 $(u,\eta )$. Moreover, the solutions we constructed
 have the following representations:
 \bq
  \label{14}
 u (z)&=&\ln \rho^I _{\vare, a^* _{ \vare} } (z)+ \vare ^{2}
  w_1 (\vare |z|) +\vare ^{2} v^* _{1,\vare} (\vare z), \\
  \label{15}
\eta (z)&=&\ln \rho^{II} _{\vare, a^* _{\vare} } (z)+ \vare ^{2}
w_2 (\vare |z|) +\vare ^{2} v^* _{2, \vare} (\vare z),
 \eq
 where the functions $\rho^{I}_{\vare , a} (z), \rho^{II}_{\vare , a}
(z)$ are defined by
 \bb
  \label{16}
  \rho^{I}_{\varepsilon ,
a}(z) = \frac{8\varepsilon^{2{N+2}} \vert f (z) \vert^2}{\lt\left(
1+\varepsilon^{2N +2} \vert F (z) +
\frac{a}{\varepsilon^{N+1}}\vert^2\right)^2} ,
 \ee
 and
 \bb
  \label{17}
  \rho^{II}_{\varepsilon ,
a}(z) = \frac{ c_0 \vare^4}{\left( 1+\varepsilon^{2N +2} \vert F
(z) + \frac{a}{\varepsilon^{N+1}}\vert^2\right)^{\lfrr}}
 \ee
 with
  \bb
 \label{18}
  f(z) = (N +1)
\prod\limits_{j=1}^{N}(z-z_{j}), \quad F (z) = \int_0^z f (\xi )
d\xi
 \ee
  for $k=1,2$,  $\vare > 0 $ and $a=a_1+i a_2\in \Bbb C$.
  The smooth radial functions, $w_1, w_2$ in (\ref{14}) and (\ref{15}) respectively
   satisfy the asymptotic formula,
  \bb
  \label{19}
  w_1 (|z|)=-C_1 \ln |z| + O(1), \qquad w_2 (|z|)=-C_2\ln |z| +O(1)
  \ee
  as $|z|\to \infty$, where
  \bq
  \label{19a}
  C_1&=&\frac{c_0\lo\lt\lf}{2(N+1)(\lt+\lf)(\lt +2\lf )},\\
  \label{19b}
  C_2&=
  &\frac{C_1\lf}{\lt}-\frac{(\lo\lf -\lt\ltt)c_0}{2(N+1)\lt }
  B\left(\frac{1}{N+1}, \lfrr -\frac{1}{N+1}\right)
  \eq
  with the beta function(Euler's integral of the first kind)
  defined by
  $$
  B(x,y)=\int_0 ^1 t^{x-1} (1-t)^{y-1} dt. \quad \forall x,y >0
  $$
  (see \cite{gra}.)
  The function $v^* _{1,\vare}, v^* _{2,\vare} $ in (\ref{14}) and
  (\ref{15}) respectively
 satisfy
 \bb
   \label{110}
 \sup_{z\in \Bbb R^2 } \frac{ |v^* _{1, \vare} (\vare z)|+|v^* _{2, \vare} (\vare z)|}{ \ln (e+|z|
 )}
  \leq o(1) \qquad \mbox{as $\vare \to 0$}.
  \ee
 \end{theorem}
 { \textsf{Remark} 1.1.} In
 the physical model of the cosmic strings of $W-$boson we note that
the coefficients in (\ref{13a}) satisfy (\ref{condition}), and  the
term containing Euler's integral
 vanishes in (\ref{19b}) to yield
$$C_2=\frac{C_1\lf}{\lt}=\frac{c_0\lo\lf ^2}{2(N+1)(\lt+\lf)(\lt +2\lf
)} >0$$
 as well as $C_1 >0$. Thus, we have extra(additional) contributions from the second terms
 of to the
 decays of $u$ and $\eta $ in (\ref{14}) and  (\ref{15})
 respectively.\\
 \ \\
 { \textsf{Remark} 1.2.} In our cosmic strings of $W-$boson we do not need smallness condition
 of the constant $G$ for the existence of condition, contrary to the other matter models of
 cosmic strings(see \cite{yan1,yan2, cha2}.)

\section{Proof of Theorem 1.1}
  \setcounter{equation}{0}
We note that for any $\vare > 0 $ and $ a \in \Bbb C $, $\ln
\rho^{I}_{\varepsilon , a}(z)$,  is a solution of the Liouville
equation(\cite{lio}).
 \bb
 \label{21}
\Delta \ln \rho^{I}_{\varepsilon , a}(z)=-\lt\rho^{I}_{\varepsilon
, a }(z)+4\pi \sum_{j=1} ^{N} \delta
 (z-z_{1,j} ).
 \ee
 We consider the following equation for $\rt $
 \bb
 \label{21a}
 \Delta \ln \rt =-\lf \ro.
 \ee
 From (\ref{21}) we have
 \bb
 \label{21b}
 \Delta \left[ \ln \ro -\sum_{j=1}^N \ln |z-z_j|^2 \right]=-\lt \ro
 .
 \ee
 Combining (\ref{21a}) with (\ref{21b}), we obtain
 $$
 \Delta \left\{ \lf \left[  \ln\ro -\sum_{j=1}^N \ln |z-z_j|^2\right] -\lt \ln
 \rt  \right\}=0,
 $$
 from which we derive
 $$
 \ln \rt =\lfr \left[ \ln \ro -\sum_{j=1}^N \ln |z-z_j|^2 \right]
 +h(z),
 $$
 where $h(z)$ is a harmonic function.
 Choosing $h(z)$ as the constant,
 $$ h(z)\equiv \lfr \ln \left(\vare ^{\frac{4\lt}{\lf} -2N-2}
 \lt^{\frac{\lt}{\lf}}[8(N+1)^2]^{-1} c_0 ^{\frac{\lt}{\lf}}\right),$$
 we get the form of $\rt$ given in (\ref{17}).
We  set
 $$
g^{I}_{\varepsilon , a}(z) =\frac{1}{\vare ^2}
\rho^{I}_{\varepsilon , a }\left(\frac{z}{\vare}\right), \quad
g^{II}_{\varepsilon , a}(z) =\frac{1}{\vare ^4}
\rho^{II}_{\varepsilon , a}\left(\frac{z}{\vare}\right),
$$
 and define $\rho_1 (r)$ and $\rho_2 (r)$ by
 $$
 \rho_1(r)=\frac{8(N +1)^2r^{2N}}{ \lt (1+r^{2N +2} )^2}
 =\lim_{\vare \to 0} g^{I}_{\varepsilon , 0}(z) ,
 $$
 and
 $$
 \rho_2(r)=\frac{c_0}{ (1+r^{2N +2} )^{\lfrr}}
 =\lim_{\vare \to 0} g^{II}_{\varepsilon , 0}(z)
 $$
 respectively.
 We transform $(u, \eta )\mapsto (v_1, v_2)$ by the
 formula
 \bb
 \label{22}
 u (z)=\ln \rho^{I}_{\varepsilon ,
a}(z) +\vare^{2} w_1 (\vare |z|) +\vare ^{2} v_1 (\vare z),
 \ee
 \bb
 \label{23}
 \eta (z) =\ln \rho^{II}_{\varepsilon ,
b}(z) +\vare^{2} w_2 (\vare |z|) +\vare ^{2} v_2 (\vare z),
 \ee
 where $w_1$ and $w_2$ are the radial functions to be determined
 below.
 Then, using (\ref{21}), the system can be written as the functional equation,
 $P(v_1,v_2, a, \vare )=(0,0)$, where
 \bb
 \label{24}
 P_1 (v_1, v_2, a,\vare )=
 \Delta v_1 + \lo \gt e^{\vare^2 (w_2 +v_2 )} +
 \lt\frac{ g^{I}_{\varepsilon , a}(z)}{\vare^2}  (e^{\vare^{2} (w_1+v_1)} -1)
+\Delta w_1,
 \ee
 and
 \bb
 \label{25}
P_2 (v_1, v_2, a, \vare )=
 \Delta v_2 +\ltt g^{II}_{\varepsilon , a}(z)
 e^{\vare^{2}(w_2+v_2)} +\lf\frac{ g^{I}_{\varepsilon , a}(z)}{\vare^2}  (e^{\vare^{2} (w_1+v_1)} -1)
  +\Delta w_2.
 \ee
  Now we introduce the functions spaces introduced in \cite{cha1}.
 For $\a >0$ the Banach spaces $X_\alpha$
and $Y_\alpha$ are defined as
\[  X_\alpha =\{ u \in L_{loc}^2 ({\mathbb R}^2) \mid
    \int_{{\mathbb R}^2} (1+|x|^{2+\alpha})|u(x)|^2 dx <\infty \} \]
equipped with the norm $\| u \|^2_{X_\alpha} =
    \int_{{\mathbb R}^2} (1+|x|^{2+\alpha})|u(x)|^2 dx$, and
\[ Y_\alpha =\{ u\in W_{loc}^{2,2}({\mathbb R}^2) \mid \| \Delta u \|_{X_\alpha}^2
  +\Big\| \frac{u(x)}{1+|x|^{1+\frac{\alpha}{2}}}\Big\|_{L^2({\mathbb R}^2)}^2 < \infty \} \]
equipped with the norm $\| u \|_{Y_\alpha}^2= \| \Delta u
\|_{X_\alpha}^2
  + \big\| \frac{u(x)}{1+|x|^{1+\frac{\alpha}{2}}}
  \big\|_{L^2({\mathbb R}^2)}^2$.
We recall the following propositions proved in \cite{cha1}.
\begin{pro}
Let $Y_\alpha$ be the function space introduced above. Then we
have the followings.
\begin{enumerate}
\item[(i)] If $v\in Y_\alpha$ is a harmonic function, then $v \equiv
 constant.$
\item[(ii)] There exists a  constant $C>0$ such that for all $v \in Y_\alpha$
    \[ |v(x)| \le C\| v \|_{Y_\alpha} \ln (e +|x|),
    \qquad \forall x\in {\mathbb R}^2 . \]
\end{enumerate}
\end{pro}
\begin{pro}
Let $\a \in (0, \frac12 )$, and let us set
 \bb
 \label{26}
 L=\Delta + \rho :\y \to \x .
 \ee
 where
 $$ \rho (z)=\rho (|z|)=\frac{8(N+1)^2 |z|^{2N}}{(1+|z|^{2N+2} )^2}. $$
We have
 \bb
 \label{27}
  Ker L=\mbox{Span} \left\{ \varphi_{+}, \varphi_{-} ,
  \varphi_{0}
 \right\},
 \ee
  where we denoted
 \bb
 \label{28}
 \varphi_+ (r,\theta)= \frac{ r^{N+1} \cos
(N+1)\theta}{1+r^{2N+2}},\quad \varphi_- (r,\theta )=
\frac{r^{N+1} \sin (N+1)\theta}{1+r^{2N+2}},
 \ee
 and
 \bb
 \label{29}
 \varphi_{0}=\frac{1-r^{2N+2}}{1+r^{2N+2}}.
 \ee
 Moreover, we have
 \bb
 \label{210}
 Im L =\{ f\in X_\alpha  | \int_{\Bbb R^2} f\varphi_{\pm}
 =0\}.
 \ee
 \end{pro}
 \ \\
 Hereafter, we fix $\alpha=\frac14$, and set $X_{\frac14} =X$ and
  $Y_{\frac14}=Y$.\\
 Using Proposition 2.1 (ii), one can check easily that for $\vare >0$
 $P$ is a well defined continuous mapping from
 $B_{\vare_0}$ into
$X ^2$, where we set $B_{\vare_0}=\{ \|v_1\|^2_{Y }+\|v_2\|^2_{Y
}+|a|^2 <\vare_0\}$, for sufficiently small $\vare_0$.
 In order to extend continuously $P$ to $\vare=0$ the radial functions $w_1 (r), w_2 (r)$
  should satisfy
 \bq
 \label{211}
 &&\Delta w_1 +\lt \rho_1 w_1 +\lo \rho_2  =0\\
 \label{212}
 &&\Delta w_2 +\lf \rho_1 w_1 +\ltt \rho_2  =0
 \eq
 For the existence and asymptotic properties of $w_1$ and $w_2$ we
 have the following lemma, which is a part of Theorem 1.1.
  \begin{lemma}
 There exist radial solutions $w_1 (|z|), w_2 (|z|)$
 of (\ref{211})-(\ref{212}) belonging to $Y $, which
  satisfy the asymptotic formula in (\ref{19}),(\ref{19a}),(\ref{19b}).
\end{lemma}

\noindent{\bf Proof:}
 Let us set $f(r)=\rho_1 (r)$. Then, it is found in \cite{bar1, cha1}
 that the ordinary differential equation(with respect to $r$),
$\Delta w_1 +C_1 \rho_1 w_1 =f(r)$ has a solution $w_1 (r)\in Y$
 given by
 \bb
 \label{213}
 w_1(r) = \varphi_0 (r) \left\{\int_0 ^r \frac{\phi_{f} (s)
-\phi_{ f}(1)}{(1-s)^2} ds + \frac{\phi_{ f }(1) r}{1-r} \right\}
 \ee
 with
 $$
 \,\, \phi_{ f} (r) := \left(\frac{
1+r^{2N+2}}{1-r^{2N+2}}\right)^2 \frac{(1-r)^2}{r} \int_0^r
\varphi_{0}(t) t{ f}(t) dt,
 $$
 where $\phi_{ f} (1)$
and $w_1(1)$ are defined as limits of $\phi_{ f} (r)$ and $w_1(r)$
as $r\to 1$. From the formula (\ref{213}) we find that
 $$
 w_1 (r)=
\varphi_0 (r) \int_2^r \left(\frac{ 1+s^{2N+2}}{1-s^{2N+2
}}\right)^2 \frac{I(s)}{s} ds +\mbox{(bounded function of $r$)}
 $$
  as $r\to
\infty$, where
$$ I(s)=  \lo \int_0^s
\varphi_0(t) t\rho_2 (t) dt.
 $$
Since $\varphi_0 (r) \rightarrow -1$ as $r\to \infty$, the first
part of (\ref{19}) follows if we show
 $$
 I =I(\infty )=\lo
  \int_0^\infty \varphi_0 (r) r\rho_2  (r) dr
 =C_1.
 $$
 Changing variable $r^{2N +2}=t$, we evaluate
 \bq
 \label{214}
 I&=& \lo \intt \varphi_0
 (r)\rho_2(r)
 rdr\n\\
 &=&c_0 \lo \intt
 \left[ \frac{r^{2N}}{(1+r^{2N_2 +2})^{3+\lfrr}}
-\frac{r^{4N+2}}{(1+r^{2N_2 +2})^{3+\lfrr}}\right]r
 dr\n\\
 &=&\frac{c_0\lo }{2(N+1)}
 \left[\intt\frac{1}{ (1+t)^{3+\lfrr} }dt-\intt \frac{t}{(1+t)^{3+\lfrr} }dt\right]\n \\
 &=&\frac{c_0\lo }{2(N+1)}\left[ \frac{1}{2+\lfrr} - \frac{1}{\left(2+\lfrr \right)\left(1+\lfrr
 \right)}\right]\n \\
 &=&
\frac{c_0\lo\lt\lf}{2(N+1)(\lt+\lf)(\lt +2\lf )} =C_1.
 \eq
 In order to obtain $C_2$ we find from (\ref{211}) and (\ref{212})
 that
 $$
 \Delta (\lf w_1 -\lt w_2 ) =(-\lo\lf +\lt \ltt )\rrt,
 $$
 from which we have
 \bq
 \label{214a}
 w_2 (z)&=&\lfr w_1 (z) +\frac{\lo\lf -\lt\ltt}{2\pi \lt}\intr \ln (|z-y|)
 \rrt (|y|)dy \n \\
 &=& -\frac{\lf C_1}{\lt} \ln|z| +\frac{\lo\lf -\lt\ltt}{2\pi \lt}\left[\intr \rrt
 (|y|)dy\right]
 \ln|z| +O(1)\n \\
 \eq
 as $|z|\to \infty$. In the case $\lo\lf -\lt\ltt=0$, we have $C_2
 =\frac{\lf C_1}{\lt}$.
 In the case $\lfrr >\frac{1}{N+1}$, we compute the integral as follows.
 \bq
 \label{214b}
 \intr\rrt (|y|)dy
 &=&2\pi c_0 \intt \frac{r}{(1+r^{2N +2})^{\lfrr}}dr\n \\
 &=& \frac{\pi c_0}{N+1}\intt \frac{t^{-\frac{N}{N+1}}}{(1+t)^{\lfrr}}
 dt \quad (r^{2N+2}=t )\n \\
 &=& \frac{\pi c_0}{N+1}B\left(\frac{1}{N+1}, \lfrr -\frac{1}{N+1}
 \right),
 \eq
 where we used the formula(See  pp.
  322\cite{gra}) for the beta function
  $$
  \intt \frac{x^{\mu -1}}{(1+x)^{\nu}} dx =B (\mu , \nu- \mu ), \qquad
   \mbox{where $\nu > \mu$}.
  $$
 Substituting (\ref{214b}) into (\ref{214a}), we have
 $w_2 (z)=-C_2 \ln |z| +O(1)$ as $|z|\to \infty$, where $C_2$ is
 given by (\ref{19b}).
This completes the proof of Lemma 2.1 $\square$\\
 \ \\

Now we compute the linearized operator of $P$.

 By direct computation we have
 \[
 \lim_{\vare \to 0} \left.\frac{\partial \go }{\partial
 a_1}\right|_{a= 0}
 =-4 \rho_1 \varphi_+ , \quad
\lim_{\vare \to 0} \left. \frac{\partial \go} {\partial
a_2}\right|_{a =0}
 =-4 \rho_1 \varphi_- ,
 \]
 \[
\lim_{\vare \to 0} \left. \frac{\partial \gt}{\partial
a_1}\right|_{a =0}
 =-4\rho_2 \varphi_+ ,\quad
\lim_{\vare \to 0}  \left.\frac{\partial \gt}{\partial
a_2}\right|_{a=0}
 =-4\rho_2 \varphi_- .
\]
Let us set $P'_{u,\eta, a }(0,0,0,0)=\mathcal{A}$. Then, using the
above preliminary computations, we obtain
 $$
 \mathcal{A}_1[\nu_1,\nu_2, \alpha ]=
  \Delta \nu_1 +\lt \rro \nu_1-4 (\lt w_1\rro +\lo \rrt )(\varphi_+ \a_1
   +\varphi_- \a_2),
  $$
   and
   $$
  \mathcal{A}_2[\nu_1,\nu_2, \alpha ]=
   \Delta \nu_2 +\lf\rro \nu_1 -4 (\lf w_1\rro +\ltt \rrt )(\varphi_+ \a_1
   +\varphi_- \a_2).
 $$
 We establish the following lemma for the operator $\mathcal{A}$.
 \begin{lemma}
The operator $\mathcal{A}:Y^2 \times \Bbb C \times \Bbb R_+ $
defined above is onto.  Moreover, kernel of $\mathcal{A}$ is given
by
 $$
 Ker \mathcal{A}=  Span\{ (0,1);
 (\varphi_\pm , \frac{\lf}{\lt} \varphi_\pm ),
 (\varphi_0 , \frac{\lf}{\lt} \varphi_0 )\}
  \times \{(0,0)\}.
 $$
 Thus, if we decompose $Y ^2\times \Bbb C=
 U
 \oplus Ker \mathcal{A}$, where we set $U=(Ker \mathcal{A})^\bot$, then
 $\mathcal{A}$ is an isomorphism from $U$ onto $X ^2$.
 \end{lemma}
 In order to prove the above lemma we need to establish the following.
 \ \\
\begin{pro}
\bb \label{214c}
 I_\pm := \intr ( \lt w_1 \rro +\lo\rrt )\varphi_\pm dx \neq 0.
 \ee
 \end{pro}
 {\bf Proof:}
 In order to transform the integrals we use the formula
 $$ L \left[
\frac{1}{16(1+r^{2N+2})^2}\right] = \frac{(N+1)^2
r^{4N+2}}{(1+r^{2N+2})^4}, \qquad \forall N\in \Bbb Z_+$$ which
can be verified by an elementary computation.  Using this,  we
have the following
 \bqn
 \lefteqn{\int_{{\Bbb R}^2} (\lt w_1 \rro +\lo\rrt)
\varphi_\pm ^2 dx = \int_0^{2\pi} \int_0 ^\infty (\lt w_1 \rro
+\lo\rrt ) \frac{r^{2N+2}}{ (1+r^{2N+2})^2}\left\{
\begin{array}{c}
 \cos^2 (N+1)\theta \\
 \sin^2 (N+1)\theta \end{array}
\right\}
rdrd\theta }\hspace{.0in}\\
  &&=\pi\int_0^\infty \left[ \frac{8(N+1)^2
r^{2N}}{(1+r^{2N+2})^2} w_1 +\lo\rrt \right]
\frac{r^{2N+2}}{(1+r^{2N+2})^2} rdr\\
&&= \pi \int_0^\infty \left[ \frac{1}{2} L  \left\{
\frac{1}{(1+r^{2N+2})^2}\right\} w_1  +\frac{\lo\rrt
r^{2N+2}}{(1+r^{2N+2})^2} \right]rdr\\
&&= \pi\int_0^\infty \left[\frac{1}{2} L w_1  \cdot
\frac{1}{(1+r^{2N+2})^2} +\frac{\lo\rrt
r^{2N+2}}{(1+r^{2N+2})^2} \right] rdr\\
&&= \pi\lo c_0
\int_0^\infty\left[-\frac{\rrt}{2(1+r^{2N+2})^2}+\frac{\rrt
r^{2N+2}}{(1+r^{2N+2})^2}\right] rdr\\
&&=\frac{\pi \lo c_0}{2}\intt \frac{r^{2N+2}
-1}{(1+r^{2N+2})^{2+\lfrr}} rdr =\frac{\pi \lo c_0}{4}\intt \frac{
t^{N +1}-1}{(1+t^{N+1} )^{2+\lfrr}}dt \quad
(\mbox{$r^2=t$})\\
&&=\frac{\pi \lo c_0}{4}\left[\int_0 ^1 \frac{ t^{N
+1}-1}{(1+t^{N+1} )^{2+\lfrr}}dt +\int_1 ^\infty \frac{ t^{N
+1}-1}{(1+t^{N+1}
)^{2+\lfrr}}dt\right]\\
&&\quad (\mbox{Changing variable $t\to 1/t$ in the second
integral,}) \\
&&=\frac{\pi \lo c_0}{4}\left[\int_0 ^1 \frac{ t^{N
+1}-1}{(1+t^{N+1} )^{2+\lfrr}}dt +\int_0 ^1
\frac{(1-t^{N+1})t^{\lfrr} }{(1+t^{N+1}
)^{2+\lfrr}}dt\right]\\
&&=\frac{\pi \lo c_0}{4}\intt \frac{(t^{N+1} -1)(1-t^{\lfrr}
)}{(1+t^{N+1} )^{2+\lfrr}}dt <0 .
 \eqn
This completes the proof of the proposition.$\square$\\

 \ \\
 We are now ready to prove Lemma 2.2.\\

\noindent{\bf Proof of Lemma 2.2:} Given $(f_1 ,f_2 )\in X ^2$, we
want first to show that there exists
 $(\nu_1, \nu_2 )\in Y ^2 $, $\a_1,\a_2\in {\Bbb R}$ such that
 $$
 \mathcal{A}(\nu_1, \nu_2, \a_1,\a_2 )=(f_1, f_2 ),
 $$
 which can be rewritten as
 \bb
 \label{214d}
  \Delta \nu_1 +\lt \rro \nu_1 -4 (\lt w_1\rro +\lo \rrt )(\varphi_+ \a_1
   +\varphi_- \a_2) =f_1,
 \ee
 and
 \bb
 \label{214e}
 \Delta \nu_2 +\lf\rro \nu_1 -4 (\lf w_1\rro +\ltt \rrt )(\varphi_+
\a_1
   +\varphi_- \a_2) =f_2 .
 \ee
 Let us set
 \bb
  \a_1 =
\frac{1}{4I_+} \int_{\Bbb R^2} f_1 \varphi_+ dx, \qquad \a_2 =
\frac{1}{4I_-} \int_{\Bbb R^2} f_2 \varphi_- dx,
  \ee
 where $I_\pm \neq 0$ is defined in (\ref{214c}). We introduce $\tilde{f}$
 by
 \bb
 \tilde{f_1}= f_1-\a_1 \varphi_+ -\a_2 \varphi_- .
 \ee
 Using  the fact
 \bq
 \int _0 ^{2\pi} \varphi_+\varphi_- d\theta
 =0,
 \eq
  we find easily
 \bb
 \int_{\Bbb R^2} \tilde{f_1}\varphi_\pm dx =0.
 \ee
Hence, by  (\ref{210}) there exists $\nu_1 \in Y$ such that
$\Delta \nu_1 +\lt \rho_1 \nu_1=\tilde{f_1}$. Thus we have found
$(\nu_1 , \a_1, \a_2)\in Y \times \Bbb R^2$ satisfying
(\ref{214d}). Given such $(\nu_1 , \a_1, \a_2)$, the function
 \bq
 \nu_2 (z)=\frac{1}{2\pi } \int_{\Bbb R^2} \ln (|z-y|) g (y) dy
 +c_1,
 \eq
 where
 $$
 g=f_2-\lf\rro \nu_1 +4 (\lf w_1\rro +\ltt \rrt )(\varphi_+
\a_1
   +\varphi_- \a_2),
   $$
   and $c_1$ is any constant,
 satisfies (\ref{214e}), and belongs to $Y$.
 We have just finished the proof that $\mathcal{A}: Y ^2 \times \Bbb
 R^2 \to X ^2$ is onto.\\
We now show that the restricted operator(denoted by the same
symbol),
 $\mathcal{A}: (Ker L\oplus
 Span\{1\})^\bot  \times \Bbb R^2 \to X ^2$ is one
 to one.
 Given $ (\nu_1,\nu_2, \a_1, \a_2)\in (Ker L\oplus
 Span\{1\})^\bot \times \Bbb R^2 $,
 let us consider the equation, $\mathcal{A}(\nu_1,\nu_2, \a_1 ,\a_2 )=(0, 0 )$, which
 corresponds to
 \bb
 \label{214f}
  \Delta \nu_1 +\lt \rro \nu_1 -4 (\lt w_1\rro +\lo \rrt )(\varphi_+ \a_1
   +\varphi_- \a_2) =0,
 \ee
 and
 \bb
 \label{214g}
 \Delta \nu_2 +\lf\rro \nu_1 -4 (\lf w_1\rro +\ltt \rrt )(\varphi_+
\a_1
   +\varphi_- \a_2) =0 .
 \ee
 Taking $L^2 (\Bbb R^2 )$ inner product of (\ref{214f}) with $\varphi_\pm$,
 and using (\ref{214c}), we find
 $\a_1=\a_2 =0$. Thus, (\ref{214f}) implies
 $\nu_1 \in Ker L$. This, combined with the hypothesis
 $\nu_1 \in (Ker L)^\bot$ leads to $\nu_1 =0$. Now, (\ref{214g}) is reduced  to
 $\Delta \nu_2=0$. Since $\nu_2\in Y$, Proposition 2.1 implies
 $\nu_2=$ constant. Since $\nu_2\in (Span\{1\})^\bot $ by hypothesis, we
 have $\nu_2=0$.
This completes the proof of the lemma.
$\square$\\
\ \\
We are now ready to prove our main theorem.\\
 \noindent{\bf Proof of
Theorem 1.1:} Let us set
\[ U =
(\mathrm{Ker} L \oplus
 \mathrm{Span}\{1\} )^\bot
 \times \Bbb \Bbb R^2 .
 \]
 Then, Lemma 2.2 shows that
 $P'_{(v_1 , v_2 , \a )} (0,0,0,0) : U \to X ^2$ is an
isomorphism. Then, the standard implicit function theorem(See e.g.
\cite{zei}), applied to the functional $P : U \times
(-\varepsilon_0, \varepsilon_0) \to X ^2$, implies that  there
exists a constant $\varepsilon_1\in (0,\varepsilon_0)$ and a
continuous function $\varepsilon \mapsto \psi ^*_\vare :=
(v_{1,\varepsilon}^*, v_{2,\varepsilon}^*, a_\varepsilon^* )$ from
$(0, \varepsilon_1)$ into a neighborhood of 0 in $U$ such that
\[ P(v_{1,\varepsilon}^*, v_{2,\varepsilon}^*, a_\varepsilon^* )
=(0,0), \quad\mbox{for all } \varepsilon\in (0, \varepsilon_1). \]
This completes the proof of Theorem 1.1.
 The
representation of solutions $u_1,u_2$, and the explicit form of
$\rho^I _{\vare, a^ * _{ \vare} } (z),$
 $ \rho^{II} _{\vare, a ^ * _{\vare} } (z),$  ,
 together with the asymptotic behaviors of $w_1, w_2$
 described in Lemma 2.1, and the fact that
 $v_{1,\varepsilon}^*, v_{2,\varepsilon}^* \in
 Y$, combined with Proposition 2.1, implies that the solutions
  satisfy the boundary condition in (\ref{13}).
 Now,
from  Proposition 2.1 we obtain that for each $j=1,2$,
 \bb
  \label{225}
  |v^*
_{j, \vare} ( x)|\leq C \Vert v^* _{j,\vare} \Vert_{Y} (\ln^+ | x|
+1) \leq C \Vert \psi_\vare \Vert_{U} (\ln^+ | x| +1).
 \ee
 This implies then
 $$
 |v^* _{j, \vare} ( \vare x)|\leq C \Vert \psi _{\vare}
\Vert_{U}(\ln^+ | \vare x| +1)  \leq C \Vert \psi_\vare
\Vert_{U}(\ln^+ | x| +1).cxxc
 $$
 From the continuity of the function $\vare\mapsto
\psi_{\vare}$ from $(0, \vare_0 )$ into $U$ and the fact $\psi^*
_0 =0$ we have
 \bb
    \label{226}
 \|\psi_\vare \|_{U} \to 0\qquad
\mbox{ as $\vare \to 0$} .
 \ee
 The proof of (\ref{110}) follows from (\ref{225})
  combined with (\ref{226}).
This
completes the proof of Theorem 1.1$\square$\\
 $$\mbox{\bf Acknowledgements} $$
This work was supported by Korea Research Foundation Grant
KRF-2002-015-CS0003.

\end{document}